\input amstex
\def\smc{\tensmc}
\documentstyle{amsppt}
\input siamptex.sty
\overfullrule=0pt
\def\Pn1{P_{n}^{(1)} (z)}
\def\Pn0{P_{n}^{(0)} (z)}

\def\cP{{\Cal P}}
\def\cQ{{\Cal Q}}
\def\greek{\alpha, \beta, \gamma, \delta, \epsilon}
\def\gree{\alpha \beta \gamma \delta \epsilon}
\def\tphin{\, {}_{10}\phi_9}
\def\ephis{\, {}_8\phi_7}

\def\fphit{\, {}_4\phi_3}
\def\tW{{\widetilde W}}

\topmatter
\vol{0}
\no{0, pp.~000--000}
\journal{SIAM J. A{\smc PPL}.M{\smc ATH}.}
\date{0}
\copyyear{1993}
\code{0}

\title WATSON'S BASIC ANALOGUE OF RAMANUJAN'S ENTRY 40
       AND ITS GENERALIZATION\endtitle

\shorttitle{Generalizing Ramanujan's Entry 40}

\author Dharma P. Gupta\fnmark{$^{\dag}$}
\and David R. Masson
\fnmark{$^{\ddag}$}\endauthor

\address{$^{\dag}$}{Department of Mathematics,
University of Toronto, Toronto, Canada, M5S 1A1}
\address{$^{\ddag}$}{Department of Mathematics, 
University of Toronto, Toronto, Canada, M5S 1A1.
Research partially supported by 
NSERC (Canada)}

\abstract{We generalize Watson's $ q $-analogue of Ramanujan's 
Entry 40 continued
fraction by deriving solutions to a  $ {}_{10} \phi_9 $  
series contiguous relation and applying Pincherle's theorem.
Watson's result is recovered as a special terminating case,
while a limit case yields a new continued fraction
associated with an  $ \ephis $  series contiguous relation.}

\subjclass 33D15, 39A10, 40A15\endsubjclass
\keywords contiguous relations, continued fractions,
Pincherle's theorem, basic hypergeometric series.\endkeywords
\endtopmatter

\centerline{\it In honour of Dick Askey on the occasion of his
60th birthday.}

\heading{1}{Introduction}
Contiguous relations for hypergeometric functions are an 
important source for obtaining explicit results for difference
equations, continued fractions, Jacobi matrices and their
corresponding orthogonal polynomials.  At the top of the 
Askey-Wilson chart of classical orthogonal polynomials [1]
one has the  $ {}_4 F_3 $  Wilson polynomials.  However the
$ {}_4 F_3 $  label is misleading, since the 
properties of these polynomials and their associated case are 
revealed by two contiguous relations for very well poised 
$ {}_7 F_6 $  series [8], [14].  These in turn can  be derived
as limits of a contiguous relation for a terminating,
very well poised, two balanced  $ {}_9F_8 $  series [18], [8].  
This  $ {}_9 F_8 $  contiguous relation is thus fundamental for the
classical hypergeometric polynomials.  In a previous publication
[15] it was shown how this  $ {}_9 F_8 $  contiguous relation
was also related to Ramanujan's famous Entry 40 continued fraction
[16],  [3].

All of the above  are   $ q \to 1 $  limits of basic hypergeometric
analogues.  Thus the  $ {}_4 \phi_3 $  Askey-Wilson polynomials
should be viewed in the light of very well poised 
$ {}_8 \phi_7 $ series [9] which are limits of terminating,
very well poised, balanced  $ {}_{10} \phi_9 $'s.  The
analogous contiguous relation for $ {}_{10} \phi_9 $'s is thus
fundamental to the whole scheme of classical and basic hypergeometric
orthogonal polynomials.  In this paper we derive this important contiguous
relation and a corresponding continued fraction.  A special
terminating version of this continued fraction yields the
following result of Watson [17] which is the  $ q $-analogue
of Ramanujan's Entry 40 [16], [3].

\thm{Theorem A {\rm (Watson [17])}} Denoting the base by
$ q^2 $ (instead of more usual  $ q $), let
$$
{1\over G(x)} \ = \ 
\prod\limits_{m=0}^\infty \ ( 1 - xq^{2m+1}), \ 
( |q| < 1 ) \ ,
$$
$$
\eqalign{
\cP  \ & = \ 
G(\alpha \beta \gamma \delta \epsilon ) 
G \Big( {\alpha\beta\gamma\over \delta\epsilon}\Big)
G \Big( {\alpha\beta\delta\over \gamma\epsilon}\Big)
G \Big( {\alpha\gamma\delta\over \beta\epsilon}\Big)
G \Big( {\alpha\beta\epsilon\over \gamma\delta}\Big)
G \Big( {\alpha\gamma \epsilon\over \beta\delta}\Big)
G \Big( {\alpha\delta\epsilon\over \beta\gamma}\Big)
G \Big( {\alpha\over \beta\gamma\delta\epsilon}\Big) \cr
\cQ \ & = \ 
G \Big({\alpha \beta \gamma \delta\over \epsilon} \Big) 
G \Big({\alpha \beta \gamma \epsilon\over \delta} \Big) 
G \Big({\alpha \beta \delta \epsilon\over \gamma} \Big) 
G \Big({\alpha \gamma \delta \epsilon\over \beta} \Big) 
G \Big({\alpha \beta \over \gamma\delta\epsilon } \Big) 
G \Big({\alpha \gamma \over \beta\delta\epsilon } \Big) 
G \Big({\alpha \delta \over \beta\gamma\epsilon } \Big) 
G \Big({\alpha \epsilon \over \beta\gamma\delta } \Big) \ .\cr}
$$
Then, provided that one of the numbers  
$ \beta, \gamma, \delta,\epsilon $  is of the form  
$ q^{\pm n} (n=1,2,\cdots ) $,
$$
{\cP - \cQ\over \cP + \cQ } \ = \ {A_0\over \beta_0 }\  
{{}\atop +} {\alpha_1\over \beta_1} \ {{}\atop + } 
{\alpha_2\over \beta_2}\ {{}\atop + \cdots } \ ,
$$
where
$$
\eqalign{
A_0 \ = \ & (q + q^{-1}) 
\Pi ( \alpha - \alpha^{-1})\cr
\alpha_m \ = \ & 
(q^{m+1} + q^{-m-1} )
(q^{m-1} + q^{1-m} )
\Pi ( \alpha^2 + \alpha^{-2} - q^{2m} - q^{-2m}) \ ,\cr
\beta_m \ = \ & 
(q^{2m+1} - q^{-2m-1} )
\Big\{ (q^m + q^{-m} )
(q^{m+1} + q^{-m-1} )
(\Sigma \alpha^2 + \Sigma \alpha^{-2} + 2) \cr
& \ - \Pi ( \alpha + \alpha^{-1}) - 
(q + q^{-1}) (q^m + q^{-m}) (q^{m+1} + q^{-m-1})
( q^{2m+1} + q^{-2m-1}) \Big\} \cr}
$$
the products and sums ranging over the numbers
$ \alpha, \beta, \gamma, \delta, \epsilon $.\endthm

A second special terminating version of the continued fraction
obtained here will give the basic analogue of Masson's 
Proposition 1 in [15], which is described as a `missing companion' of
Ramanujan's Entry 40.  For the sake of completeness we
state Masson's result:

\thm{Theorem B {\rm (Masson [15])}}  Let  
$ P' = \Pi \Gamma  \big(( 3+\alpha\pm \beta \pm \gamma \pm \delta
\pm \epsilon)/4\big)$
($0,2 $  or  $ 4 $  minus signs) and
$ Q' = \Pi \Gamma  \big(( 1+\alpha\pm \beta \pm \gamma \pm \delta
\pm \epsilon)/4\big)$
($1$ or $ 3 $  minus signs).
Then if one of the parameters  $ \beta , \gamma, \delta, \epsilon $
is an odd integer,
$$
{Q'\over P'} \ = \ 
 {-1\over a_0} \ {{}\atop -} {2b_1\over a_1 }\  {{}\atop -}
{b_2\over a_2} \  {{}\atop -} 
{b_3\over a_3} \  {{}\atop - \cdots } 
$$
where
$$
\eqalign{
b_n  =  &
\Big( \Pi \big( (2n-1)^2 -\alpha^2\big)\Big)\Big/
(16)^3 (2n-1)^2, \cr
a_n = &
\Big\{ 2n^6 + n^4 ( 5 - \Sigma \alpha^2 )\big/ 4+n^2
\big( - 26 + ( 1 + \Sigma \alpha^2)^2 - 2 \Sigma \alpha^4\big)\big/
64 - a_0\Big\} \Big/(4 n^2 -1)\, , \cr
a_0  =  &
\big\{ 2 ( 1 - \Sigma \alpha^4) + (1 -\Sigma \alpha^2)^2 -
8\Pi \alpha \big\} \Big/ (16)^2 \, ,\cr}
$$
with these products and sums ranging over the parameters
$ \greek $.\endthm

Masson [15] also gave the non-terminating versions of
Ramanujan's Entry 40 and Theorem B.

The object of the present study is to obtain the non-terminating
versions of Watson's theorem and the  $ q $-analogue of Masson's
theorem given above.  They are given in Section 4 by Corollaries
7 and 8 respectively. Our approach is similar to that in several
recent papers [5], [6], [12], [13] on the subject where
Pincherle's theorem [11] has been used  to bring out the connection
between  several of Ramanujan's Chapter 12 entries and the
general theory of hypergeometric orthogonal functions (Askey and
Wilson [1], Wilson [18]). For other approaches to explaining
some of Ramanujan's continued fraction entries see [3], [10], [19].

\heading{2}{Contiguous relation}
We consider a terminating very well poised  balanced
$ {}_{10}\phi_9 $  basic hypergeometric function
$$
\leqalignno{
\phi   =  & \phi ( a; b, c,d,e,f,g,h) \cr
       :=  & \tphin \left (
       { {a,q\sqrt a, -q\sqrt a, b, c,d,e,f,g,h;\qquad q,q }
       \atop {\displaystyle \sqrt a, -\sqrt a, {aq\over b}, 
       {aq\over c}, {aq\over d}, {aq\over e},
    {aq\over f}, {aq\over g}, {aq\over h} \phantom {q,q}}}\right)\ ,
    & (2.1) \cr
    a^3 q^2  =  &  bcdefgh\ , \qquad |q| < 1 \cr}
$$
with say  $ h = q^{-n} $,  $ n= 0,1, \, \ldots $, 
$ g = sq^{n-1} $,
$ s := {a^3 q^3 \over bcdef} $.
We follow the usual notation for variations of  $ \phi $  with
respect to the parameters.  For example  $ \phi ( b + , c - ) $
represents the  $ \phi $  with  $ b $  and  $ c $  replaced by
$ bq $  and  $ {c\over q} $  respectively.  $ \phi_+ $
denotes the  $ \tphin $  got by replacing  $ a $  by  $ aq^2 $  and
$ b,c,d,e,f,g,h $  by  $ bq, cq, dq , eq, fq, gq, hq $
respectively.  

We need a contiguous relation basic analogue to the contiguous
relation derived by Wilson [18] for the  $ {}_9 F_8 $
hypergeometric function.  In order to work out this contiguous
relation, we shall use Wilson's method [18] using the basic
hypergeometric analogues of the relevant formulas.

\thm{Lemma 1}  Let  $ \phi $  be given by (2.1) (not necessarily
terminating).  Then
$$
\leqalignno{
& \phi(b-, c+) - \phi & (2.2) \cr
= \ & {{aq\over c} (1- {cq\over b}) (1-{bc\over aq}) (1-aq)
(1-aq^2) (1-d)(1-e)(1-f)(1-g)(1-h)\over
(1-{aq\over b}) (1-{aq^2\over b}) ( 1-{a\over c})
(1-{aq\over c}) (1- {aq\over d}) (1- {aq\over e}) 
(1-{aq\over f})
(1-{aq\over g})
(1-{aq\over h})} \cr
& \hskip 2in \times \phi_+ (b-) \ .\cr}
$$\endthm

\prf{Proof} A straightforward term by term subtraction on the left side
of (2.2) leads to the result. \qquad\endproof

\thm{Lemma 2} If  $\phi$  (given by (2.1)) is terminating,
then
$$
\leqalignno{
& {b^2 (1-h) 
(1-{aq\over bc}) (1-{aq\over bd}) (1-{aq\over be})
(1-{aq\over bf}) (1-{aq\over bg}) \over
(1-{aq\over b}) (1-{aq^2\over b})}\ \phi_+ (b-) & (2.3)\cr
- & {h^2 (1-b)
(1-{aq\over ch}) (1-{aq\over dh}) (1-{aq\over eh})
(1-{aq\over fh}) (1-{aq\over gh}) \over
(1-{aq\over h}) (1-{aq^2\over h})}\ \phi_+ (h-) \cr
- & {b (1- {h\over b})
(1-{aq\over c}) (1-{aq\over d}) (1-{aq\over e})
(1-{aq\over f}) (1-{aq\over g}) \over
(1-aq) (1-aq^2)}\ \phi \ = \ 0\ . \cr}
$$\endthm

\prf{Proof} By eliminating  $ \phi_+ ( b-) $  from (2.2) and another
similar relation written for  $ \phi ( b- , d+ ) - \phi $
we obtain
$$
\leqalignno{
& c ( 1-c) 
( 1-{a\over c}) ( 1- {dq\over b}) (1-{bd\over aq})
\phi ( b-, c+ ) & (2.4) \cr
-\ & d ( 1 - d) 
( 1-{a\over d}) ( 1- {cq\over b}) (1-{bc\over aq})
\phi ( b-, d+ ) \cr
+\  & d ( 1 -{b\over q} ) 
( 1-{c\over d}) ( 1- {aq\over b}) (1-{cd\over a})
\phi\ =\ 0\ . \cr}
$$
With say,  $ h = q^{-n} $,  we can apply an iterate of Bailey's
transformation  8.5(1) [2, p. 68] to  $ \phi $,  $ \phi_+ (b-) $  and
$ \phi_+ ( h- ) $ (the transformation [4, exercise 2.19, p. 53] with  
$ b,e,g $  replaced by  $ g, b, e $  respectively). The three
transformed series are related via (2.4).  
Reversing the transformations in this relation we arrive at (2.3). 
\qquad\endproof

\thm{Theorem 3} If  $ \phi $  (given by (2.1)) is terminating,
then
$$
\leqalignno{
& {g(1-h)
( 1-{a\over h}) ( 1- {aq\over h}) (1-{aq\over gb})
( 1-{aq \over gc}) ( 1- {aq\over gd})
( 1-{aq\over ge}) ( 1- {aq\over gf}) 
\over
(1-{hq\over g})} \cr 
& \qquad \times  
[\phi (g-, h+) - \phi]  & (2.5) \cr
- &  {h(1-g)
( 1-{a\over g}) ( 1- {aq\over g}) (1-{aq\over hb})
( 1-{aq \over hc}) ( 1- {aq\over hd})
( 1-{aq \over he}) ( 1- {aq\over hf})
\over
(1-{gq\over h})} \cr
& \qquad \times  
[\phi (h-, g+) - \phi] \cr
- &  {aq\over h} 
( 1-{h \over g}) ( 1- {gh\over aq})
(1-b) (1-c) (1-d) (1-e) (1-f) \phi \ = \ 0 \ . \cr}
$$\endthm

\prf{Proof} We eliminate  $ \phi_+ (b-) $  and  $ \phi_+ (c-) $
from (2.2), (2.2) with  $ b \leftrightarrow c $  and (2.3)
with  $ c\leftrightarrow h $.  A final interchange of parameters
$ b \leftrightarrow g $,  $ c \leftrightarrow h $  yields
the desired result. \qquad\endproof

Substituting  $ h = q^{-n} $,  $g = sq^{n-1} $,  and
renormalizing, the contiguous relation (2.5) becomes the linear
second order difference equation
$$
 X_{n+1} - a_n X_n + b_n X_{n-1} \ = \ 0 \ , \qquad
n \ge 0  
\leqno(2.6)
$$
$$
\leqalignno{
a_n \ = \ & \Big[ {q^{-n+{1\over 2}} \over \sqrt s} 
(1- sq^{n-1})
(1- {s\over a} q^{n-1})
(1- {s\over a} q^{n-2 }) \cr
& \qquad \times
{(1- {a\over b} q^{n+1})
(1- {a\over c} q^{n+1})
(1- {a\over d} q^{n+1})
(1- {a\over e} q^{n+1})
(1- {a\over f} q^{n+1})\over
(1-s q^{2n})} \cr
& + {q^{-n+{3\over 2}}\over \sqrt s} 
(1-q^n)(1-a q^n) 
(1-aq^{n+1}) \cr
& \qquad \times
{(1- {bs\over a} q^{n-2})
(1- {cs\over a} q^{n-2})
(1- {ds\over a} q^{n-2})
(1- {es\over a} q^{n-2})
(1- {fs\over a} q^{n-2})\over 
(1-s q^{2n-2})} & (2.7) \cr
& + {\sqrt s\over a} q^{n-{1\over 2}}
(1-s q^{2n-1}) 
(1-{s\over aq^2})
(1-b)(1-c)(1-d)(1-e)(1-f) \Big] \cr
& \qquad \Big/ \big[(1-sq^{2n-1}) ( 1-{s\over a} q^{n-2} )
(1-aq^{n+1})\big] \ , \cr
b_n \ = \ & 
{q^{-2n+3} \over s} 
(1-q^n) (1-sq^{n-2} )
(1- {a\over b} q^n) 
(1- {a\over c} q^n) 
(1- {a\over d} q^n) 
(1- {a\over e} q^n) 
(1- {a\over f} q^n)  \cr
&\ \times 
{(1- {bs\over a} q^{n-2})
(1- {cs\over a} q^{n-2})
(1- {ds\over a} q^{n-2})
(1- {es\over a} q^{n-2})
(1- {fs\over a} q^{n-2})\over 
(1-s q^{2n-1})
(1-s q^{2n-2})^2
(1-s q^{2n-3})} \ , & (2.8)\cr}
$$
$$
s \ = \ {a^3 q^3 \over bcdef }
$$
with the solution
$$
\leqalignno{
X_n^{(1)} \ = \ & {q^{-{n^2\over 2} + n}\over s^{n\over 2}}\ 
{(sq^{2n-1})_\infty 
(aq^{n+1})_\infty \over
(sq^{n-1})_\infty 
({s\over a}q^{n-1})_\infty 
\big( {aq^{n+1}\over b}, 
{aq^{n+1}\over c}, 
{aq^{n+1}\over d}, 
{aq^{n+1}\over e}, 
{aq^{n+1}\over f}\big)_\infty } & (2.9)\cr
& \times \phi (a; b,c,d,e,f, sq^{n-1} , q^{-n} ) \ . \cr}
$$
Here the infinite product  $ (a)_\infty $  means
$$
\eqalignno{
& (a)_\infty \ = \ (a;q)_\infty \ = \ (1-a)(1-aq) (1-aq^2) \ \ldots \cr
\noalign{\hbox{and}} \cr
&  (a,b, \, \ldots, k)_\infty \ = \ (a)_\infty (b)_\infty \ldots
(k)_\infty \ . \cr}
$$
For the exceptional values  $ s = q,q^2 $, the  $ a_n, b_n $  and
$ b_{n+1} $  in  (2.7), (2.8) and
$ X_n^{(1)} , X_{n-1}^{(1)} $  in (2.9) are indeterminate
at  $ n = 0 $.  We resolve this indeterminancy by taking limits
as  $ n \to 0 $.

Next, we proceed to find a second linearly independent solution
to the second order difference equation (2.6).  This can be obtained
by using a  $ q$-analogue of [15].  Thus from (2.6), (2.9)
and a symmetry relation ((2.11), below) we are able to
obtain a second terminating  $ \tphin $ solution for the special
values  $ s = q, q^2 , \, \ldots $.  For general values of  $ s $
the second solution will be an appropriate combination of two 
non-terminating  $ \tphin $'s which satisfy a four-term
transformation (Gasper and Rahman [4], formula III.39,
p. 247).  We will consider the case of general  $ s $  in future
work.

Observe that with the replacement
$$
(a,b,c,d,e,f,sq^{n-1}, q^{-n} )\  \to \ 
\Big({q\over a},\ {q\over b},\ {q\over c},\ {q\over d},\ {q\over e},\ 
{q\over f},\ {q^{-n+2}\over s }, \ q^{n+1} \Big) 
\leqno(2.10)
$$
we have
$$
(a_n, b_n) \ \to \ (a_n, b_{n+1}) \ .
\leqno(2.11)
$$
It is easy to check that  $ b_n \to b_{n+1}$.  To check
$ a_n \to a_n $  we used the `Maple' software on the
computer.  This meant verifying a polynomial identity in
$ x = q^{-n} $  of degree fourteen.

Applying the transformation (2.10) to (2.6) and (2.9)
and renormalizing, we obtain the second solution
$$
\leqalignno{
 X_n^{(2)}  =  & 
{q^{-{n^2\over 2} + n} \over s^{n\over 2}}
{({s\over a}q^n)_\infty 
(sq^{2n-1})_\infty \over  
(q^{n+1})_\infty
(aq^n)_\infty 
({bs\over a} q^{n-1},
{cs\over a} q^{n-1},
{ds\over a} q^{n-1},
{es\over a} q^{n-1},
{fs\over a} q^{n-1})_\infty } \cr
& \times \phi \Big({q\over a}; 
{q\over b}, {q\over c},{q\over d},{q\over e},{q\over f},
{q^{-n+2}\over s} , q^{n+1} \Big) \ ,& (2.12) \cr
 & \hskip 1.25in s  = q , q^2  , \, \ldots \ . \cr}
 $$
 Note that  $  \phi $  is terminating in (2.12) because of the
 parameter  $q^{-n+2}/s $.

\heading{3}{Asymptotics and Pincherle's theorem}
In order to obtain a minimal (subdominant) solution for (2.6) we need
the large  $ n $  asymptotics of (2.9) and (2.12).  Applying
Tannery's theorem to the ${}_{10} \phi_9 $'s on the right side
of (2.9) and (2.12)  we have, as  $ n \to \infty $,
$$
\leqalignno{
& X_n^{(1)} \sim
{q^{-{n^2\over 2} + n} \over s^{n\over 2}}\ 
W (a; b,c,d,e,f) \ , \qquad 
\Big| {s\over aq } \Big| < 1 \ , & (3.1)\cr
\noalign{\hbox{and}} \cr
& X_n^{(2)} \sim
{q^{-{n^2\over 2} + n} \over s^{n\over 2}}\ 
W \Big({q\over a}; {q\over b},{q\over c},{q\over d},{q\over e},
{q\over f}\Big) \ , \qquad 
\Big| {aq^2\over s } \Big| < 1 \ , & (3.2)\cr}
$$
where
$$
W (a; b,c,d,e,f) \ := \ 
\ephis 
\left (  
       { {a,q\sqrt a, -q\sqrt a, b, c,d,e,f} 
       \atop 
       {\displaystyle \sqrt a, -\sqrt a, {qa\over b}, {qa\over c},
       {qa\over d}, {qa\over e},
    {qa\over f}}} \ ; \  
    q, {a^2 q^2\over bcdef} \right)  \  .
$$
We write
$$
W_1 \ := \ W (a; b,c,d,e,f) \ , \qquad 
|s| < | qa| 
$$
and its analytic continuation otherwise; and
$$
W _2 \ :=\  W
\Big({q\over a}; {q\over b},{q\over c},{q\over d},{q\over e},
{q\over f}\Big) \ , \qquad 
|s| > | aq^2| 
$$
and its analytic continuation otherwise.

Taking now
$$
X_n^{(3)} \ := \ W_2 X_n^{(1)} - W_1 X_n ^{(2)}
\leqno(3.3)
$$
it follows from (3.1) and (3.2) that
$$
\lim\limits_{n\to \infty} \ 
{X_n^{(3)}\over X_n^{(1)}} \ = \ 0 \ ,
\Big| {s\over q}\Big| < |a| < \Big|{s\over q^2} \Big| \ .
\leqno(3.4)
$$
This establishes that  $ X_n^{(3)} $  is a minimal solution of
(2.6).  An application of Pincherle's theorem [11] then leads to 
the following result:

\thm{Theorem 4}  Let  $ s = q,q^2 , \, \ldots $.  Then
$$
{1\over a_0} \ {{}\atop -} 
{b_1\over a_1} \ {{}\atop -} 
{b_2\over a_2} \ {{}\atop - \cdots }  \ = \ \lim\limits_{n\to 0} \ 
{W_2 
X_n^{(1)} - W_1 X_n^{(2)} \over
b_n ( W_2 X_{n-1}^{(1)} - W_1 X_{n-1}^{(2)})} \ .
\leqno(3.5)
$$\endthm

\prf{Proof} From Pincherle's theorem (3.5)  is true for 
$ \big|{s\over q}\big | < |a| < \big| {s\over q^2} \big| $.
For other values of  $ a $  the result follows
by analytic continuation.  To the left side of (3.5)
we can apply the `parabola theorem' (see Jones and Thron 
[11, p. 99] and Jacobsen [10]), since from (2.6),
$ {b_n\over a_n a_{n-1}} = {q^3\over (1+q)^2} \big( 1 + O(q^n)\big) $.
Hence the left side of (3.5) is a meromorphic function of $a$.
The right side of (3.5) involves convergent infinite products and
$ \ephis $'s  which are each expressible in terms of convergent
infinite products and convergent   $ \fphit $'s (Gasper and 
Rahman [4, (2.10.10), page 43]).
Consequently the right side of (3.5) is also a meromorphic function
of $a$ and (3.5) follows by analytic continuation to all values
of $a$.
Note that the exceptional cases  
$ s = q, q^2 $  
which
cause indeterminancy are taken care of by the limit
$ n\to 0 $  on the right
side of (3.5). \qquad\endproof

For the exceptional values
$ s = q^2, q $,
the above theorem gives respectively the non-terminating versions
of Theorem A (Watson[17]) and the basic analogue of Theorem B
(Masson [15]).  We now demonstrate how to derive the
terminating versions of Theorem 4.  We shall need to express the
ratio  $ W_1/W_2 $  in terms of infinite products when
$ {b\over a} = q^N $,  $ N $  being an integer.  We write
$$
\tW (a; b,c,d,e,f) \ := \ 
\left ( {aq\over b}, {aq\over c}, {aq\over d}, {aq\over e},
    {aq\over f}\right)_\infty   
W (a; b,c,d,e,f)
$$
and
$$
U (a; b,c,d,e,f) \ := \ 
{\tW (a; b,c,d,e,f)\over (aq, b,c,d,e,f)_\infty }\ .
$$

\thm{Lemma 5}  If  $ {b\over a} = q^N $,  where  $N $  is an
integer, then
$$
U (a; b,c,d,e,f) \ = \ 
\big({s\over aq}\big)^N 
U \Big({b^2\over a} ; b, {bc\over a}, {bd\over a} ,
{be\over a} ,{bf\over a}\Big) \ .
\leqno(3.6)
$$\endthm

\prf{Proof} Refer to Bailey's three-term  $ \ephis $  transformation
(formula III.37, p. 246 of Gasper \& Rahman [4]).  If we apply
the condition   $ {b\over a} = q^N $,  $ N= 0 $,  $ \pm 1 $,
$\pm 2 , \ldots $ we obtain the desired result.\qquad\endproof

\thm{Lemma 6}  If  
$ s = {a^3 q^3\over bcdef} = q^M $  and  $ {b\over a} = q^N $,  
$ M $  and  $ N $  being integers, then
$$
\leqalignno{
&{\tW (a; b,c,d,e,f)\over \tW \big(
 {q\over a}; {q\over b}, {q\over c}, {q\over d}, {q\over e},
    {q\over f}\big)} & (3.7) \cr
= \ &
\lambda {(aq, c,d,e,f, 
{aq^2\over s},
{aq\over ef}, 
{aq\over df}, 
{aq\over de}) _\infty \over
( {bc\over a}, {bd\over a}, {be\over a}, {bf\over a},
{q^2\over a} ,  {q\over b}, {cd\over a}, {ce\over a}, 
{cf\over a})_\infty } \ , \cr}
$$
where
$$
\leqalignno{
&\qquad \lambda  = (-1)^{n+1} \big( {s\over aq}\big)^N
\big( {c\over b}\big)^{n+1} q^{n(n+1)/2} \quad \hbox{for} \quad
{aq^3\over bs} = q^{-n}, \ n=0,1,2,\, \ldots , & (3.8)\cr
\noalign{\hbox{and}} \cr
&\qquad \lambda  = (-1)^{n+1} \big( {s\over aq}\big)^N
\big( {b\over c}\big)^{n+1} q^{(n+1)(n+2)/2} \quad \hbox{for} \quad
{bs\over aq} = q^{-n}, \ n= -1,0,1,2,\, \ldots \, . & (3.9)\cr}
$$\endthm

\prf{Proof} We express the left side of (3.7) in terms of appropriate
$ {}_4 \phi_3 $'s.
To the numerator  $ W $  in (3.7) we first apply the identity (3.6)
and then the three-term transformation formula III.36
[4, p. 246].  To the denominator  $ W $  we first apply  the
$ \ephis $  transformation formula III.24 [4, p. 243] and then
the formula III.36.  We also make use of the relation
$$
\leqalignno{
\lim\limits_{e\to q^{-n}} \ (e)_\infty 
\fphit 
\Big({{a,b,c,d}\atop {e,f,g}}; \ q,q \Big) \ = \ 
& {(a,b,c,d, fq^{n+1}, gq^{n+1} , q^{n+2})_\infty \over
(aq^{n+1}, bq^{n+1}, cq^{n+1}, dq^{n+1}, f,g)_\infty }\cr
& \quad \times q^{n+1} \fphit
\Big({{aq^{n+1},bq^{n+1},cq^{n+1},dq^{n+1}}\atop 
{q^{n+2}, fq^{n+1},gq^{n+1}}}; \ q,q \Big) \ .& (3.10)\cr}
$$
All this enables us to recognize and cancel a common linear combination of
$ \fphit$'s from the numerator and the denominator yielding
the desired result.
We note that in the case  $ {bs\over aq } = q$, $ s= q $,  the
limit (3.10) is not required and there is an exact cancellation.  
\qquad\endproof

\heading{4}{Exceptional values  $ s = q,q^2 $}
We now restate Theorem 4 for the exceptional values 
$ s = q,q^2 $  and the form they take when the continued
fraction terminates:

\thm{Corollary 7}  If  $ s = q^2 $,  then (3.5) can be 
rewritten as
$$
{1\over a_0} \ {{}\atop -}
{b_1\over a_1} \ {{}\atop -} 
{b_2\over a_2} \ {{}\atop - \cdots}  \ = \ 
{2a(1-q) \over q^{3/2} (1-\alpha) (1-\beta)(1-\gamma)(1-\delta)
(1-\epsilon)} \ 
\Big( {1-V\over 1+V } \Big) \ , 
\leqno(4.1)
$$
$$\leqalignno{ 
a_n \ = \  & \Big[ q^{1\over 2} \Pi 
(\alpha^{1\over 2} + \alpha^{-{1\over 2}}) +
q^{1\over 2} 
(q^{1\over 2} + q^{-{1\over 2}} )
(q^{n+{1\over 2}} + q^{-n-{1\over 2}}) \cr
&  \times
(q^{n\over 2} + q^{-{n\over 2}} )
(q^{n+1\over 2} + q^{-{n+1\over 2}})  & (4.2)\cr
& \quad - q^{1\over 2} 
(\Sigma \alpha + \Sigma \alpha^{-1} + 2) 
(q^{n\over 2} + q^{-{n\over 2}} )
(q^{n+1\over 2} + q^{-{n+1\over 2}}) \Big] \cr
& \qquad \Big/
(q^{n\over 2} + q^{-{n\over 2}} )
(q^{n+1\over 2} + q^{-{n+1\over 2}})\ , \cr
b_n \ = \ &
{- q\Pi (\alpha + \alpha^{-1} - q^n - q^{-n})\over
(q^{-{n\over 2}} + q^{{n\over 2}} )^2
(q^{-n-{1\over 2}} - q^{n+{1\over 2}})
(q^{-n+{1\over 2}} - q^{n-{1\over 2}})} \ , \cr
V \ = \ & { ({q\over a})_\infty 
 ({q^2\over a})_\infty \over (a)_\infty (aq)_\infty } \ 
 {\tW_1\over \tW_2 } \ , & (4.3) \cr
a \ = \ & ( q\gree )^{1\over 2} \ , \cr}
$$
and product  $ \Pi $  and summation  $ \Sigma $  are taken
over the parameters  $ \greek $.  If one of the parameters
$ \beta, \gamma , \delta , \epsilon = q^N $,  $ N $  integer,
$ ( \greek ) \to ( \alpha^2 , \beta^2, \gamma^2 , \delta^2,
\epsilon^2) $  and the base  $ q $  is changed to
$ q^2 $,  then the right side of (4.1) becomes
$$
{2 (q^{-1} - q) \over q	\Pi (\alpha - \alpha^{-1} )} \ 
{\cP - \cQ\over \cP + \cQ }
\leqno(4.4)
$$
with
$$
\leqalignno{
{1\over \cP } \ = \ 
& \Big(\gree , 
{\alpha\over \beta \gamma \delta \epsilon},
{\alpha \beta \gamma \over \delta \epsilon},
{\alpha \beta \delta \over \gamma \epsilon},
{\alpha \gamma \delta \over \beta \epsilon},
{\alpha \beta \epsilon \over \gamma \delta},
{\alpha \gamma \epsilon \over \beta \delta},
{\alpha \delta \epsilon \over \beta \gamma};\ q^2 \Big)_\infty  \cr
{1\over \cQ } \ = \ 
& \Big( {\alpha \beta \gamma \delta\over \epsilon } ,
{\alpha \beta \gamma \epsilon\over \delta },
{\alpha \beta \delta \epsilon\over \gamma },
{\alpha \gamma \delta \epsilon\over \beta },
{\alpha \beta\over \gamma \delta \epsilon},
{\alpha \gamma\over \beta \delta \epsilon},
{\alpha\delta\over \beta\gamma \epsilon } ,
{\alpha\epsilon\over \beta\gamma \delta } ;\ q^2 \Big)_\infty \cr}
$$
and $ a_n , b_n $  modified accordingly.  This reproduces Theorem A.
\endthm

\prf{Proof} We write  $ \alpha = {a\over b}$,
$ \beta  = {a\over c}$,
$ \gamma  = {a\over d}$,
$ \delta  = {a\over e}$,
$ \epsilon  = {a\over f}$.
After that the proof is straightforward on substituting the values
of  
$ X_0^{(1)} $,  $ X_0^{(2)} $, 
$ b_0 X_{-1}^{(1)} $,  $ b_0 X_{-1}^{(2)} $ 
from (2.9), (2.12) and (2.8) into Theorem 4.  A lot
of algebra is involved in the simplification.  Also, appropriate
limits are to be taken whenever indeterminates occur.

In order to obtain the terminating form (4.4) we need to use 
Lemma 6, after interchanging say  $ b $  and  $ f $.  In both
cases viz., 
$$ 
{aq^3 \over fs}\  = \ q^{-n} \ ,  \quad n=0,1,2, \ldots \qquad
\hbox{and} \qquad  
{fs\over aq}\  = \ q^{-n}\  , \quad n = -1, 0,1,2 \, \ldots 
$$
whether the termination is due to one or the other, the result
works out to be the same.  The above result (4.4) yields 
Watson's result [17] i.e. Theorem A in Section 1.\qquad\endproof

\thm{Corollary 8}  If  $ s = q $,  then (3.5) can be 
rewritten as
$$
{1\over a_0} \ {{}\atop -}
{b_1\over a_1} \ {{}\atop -} 
{b_2\over a_2} \ {{}\atop - \cdots}  \ = \ 
2\left(
a_0 + {a^2\over q} \ 
{({1\over a})_\infty ({1\over a})_\infty\over
 (aq)_\infty ({a\over q})_\infty } \ 
{\tW_1\over\tW_2} \right)^{-1} \ , 
\leqno(4.5)
$$
$$
\leqalignno{
a_n \ = \ & 
{q^{1\over 2}\over 
(q^{-n-{1\over 2}} - q^{n+{1\over 2}}) 
(q^{-n+{1\over 2}} - q^{n-{1\over 2}}) } \ 
\bigg\{ (q^{3n} +  q^{-3n})
(q^{1\over 2} + q^{-{1\over 2}}) & (4.6)\cr
& - (q^{2n} + q^{-2n}) 
\Big[ q^{1\over 2} + q^{-{1\over 2}} + 
\Sigma (\alpha^{-{1\over 2}} + \alpha^{1\over 2})\Big] \cr
& + (q^{n} + q^{-n}) 
\Big[ - q^{3\over 2} - q^{-{3\over 2}} + 
(q^{1\over 4} + q^{-{1\over 4}} )
\Pi  (\alpha^{1\over 4} + \alpha^{-{1\over 4}})\cr
& \qquad + 
(q^{1\over 4} - q^{-{1\over 4}} )
\Pi  (\alpha^{1\over 4} - \alpha^{-{1\over 4}})\Big]\cr
& + (q^{-1} + q) 
(q^{1\over 2} + q^{-{1\over 2}}) 
 + (q^{-1} + q) 
\Sigma (\alpha^{-{1\over 2}} + \alpha^{1\over 2}) \cr
& - 
(q^{1\over 4} + q^{-{1\over 4}})
(q^{1\over 2} + q^{-{1\over 2}}) 
\Pi  (\alpha^{1\over 4} + \alpha^{-{1\over 4}})\cr
& + 
(q^{1\over 4} - q^{-{1\over 4}})
(q^{1\over 2} + q^{-{1\over 2}}) 
\Pi  (\alpha^{1\over 4} - \alpha^{-{1\over 4}})\bigg\} \ , \cr
b_n \ = \ & 
q^{5\over 2} \ 
{\Pi (q^{-n} + q^n - \alpha^{1\over 2} q^{-{1\over 2}} -
\alpha^{-{1\over 2}} q^{1\over 2}) \over
(q^{n\over 2} + q^{-{n\over 2}}) 
(q^{{n\over 2}-{1\over 2}} + q^{-{n\over 2} +{1\over 2}}) 
(q^{n-{1\over 2}} - q^{-n+{1\over 2}})^2 } \ , \cr}
$$
$$
\leqalignno{
\tW_1 \ = \ &  
\tW
\Big( a;\  a \sqrt {q\over \alpha}, \ 
a \sqrt {q\over \beta}, \ 
a \sqrt {q\over \gamma}, \ 
a \sqrt {q\over \delta}, \ 
a \sqrt {q\over \epsilon} \Big) \ , & (4.7)\cr
\tW_2 \ = \ &  
\tW
\Big({q\over a} ; \  {\sqrt {\alpha q}\over a}, \ 
 {\sqrt {\beta q}\over a }, \ 
 {\sqrt {\gamma q}\over a}, \ 
 {\sqrt {\delta q}\over a }, \ 
 {\sqrt {\epsilon q}\over a}  \Big) \ , \cr
a \ = \ &
\Big( {\gree\over q} \Big)^{1\over 4} \ , \cr}
$$
and product  $ \Pi $  and summation  $ \Sigma $  are taken over
parameters  $ \greek $.  If one of the parameters
$ \beta , \gamma, \delta, \epsilon $  is  $ q^N $,  $ N $  an odd
integer,  $ (\greek )  \to (  \alpha^4, \beta^4, \gamma^4 , \delta^4 ,
\epsilon^4 ) $  and the base $ q $  is changed
to  $ q^4 $, then (4.5) becomes
$$
{1\over a_0} \ {{}\atop -}
{b_1\over a_1} \ {{}\atop -} 
{b_2\over a_2} \ {{}\atop - \cdots}  \ = \ 
2\left(
a_0 - {q^2\over \alpha^2} \ 
{\cP'\over\cQ'} \right)^{-1} \ ,  
\leqno (4.8)
$$
where
$$
\leqalignno{
& \quad {1\over \cQ'}  =  \bigg( 
{q\alpha \beta \gamma \delta \over \epsilon } ,  
{q\alpha \gamma \delta \epsilon \over \beta } ,  
{q\alpha \beta \epsilon \delta \over \gamma } ,  
{q\alpha  \epsilon \beta \gamma \over \delta} , 
{q\alpha  \beta \over \epsilon \gamma \delta} , 
{q\alpha \gamma \over \epsilon \beta \delta} , 
{q\alpha \delta \over \epsilon \beta\gamma } , 
{q\alpha  \epsilon \over \beta \gamma \delta} ; 
q^4 \bigg)_\infty \, , & (4.9) \cr
& \quad {1\over \cP'}  =  \bigg( 
q^3 \gree , 
{q^3\alpha \over \beta \gamma \delta \epsilon } , 
{q^3\alpha \delta \epsilon \over \beta \gamma } , 
{q^3\alpha \gamma \epsilon \over \beta \delta } , 
{q^3\gamma \delta \alpha \over \epsilon \beta } , 
{q^3\alpha \beta  \delta\over  \epsilon \gamma } , 
{q^3\alpha \beta \gamma\over  \epsilon \delta } , 
{q^3\alpha \beta  \epsilon \over \gamma\delta } ; 
q^4 \bigg)_\infty \, ,  \cr}
$$
with  $ a_n  $  and  $ b_n $  modified accordingly.\endthm

\prf{Proof} We write  
$ \alpha = {a^2q\over b^2} $,
$ \beta = {a^2q\over c^2} $,
$ \gamma = {a^2q\over d^2} $,
$ \delta = {a^2q\over e^2} $ 
and  
$ \epsilon = {a^2q\over f^2 } $ 
and make the appropriate substitutions from (2.9), (2.12) and (2.8)
into Theorem 4.  A considerable amount of algebra is required to 
reexpress (2.7) as (4.6) for which we used `Maple' software on the computer.
We use Lemma 6 after interchanging say  $ b $  and  $ f $  to
arrive at (4.8).  The result is the same for either type of
termination (3.8) or (3.9).  Note that (4.8) can be reexpressed in the form
$$
- {q^2\over \alpha^2 } \ {\cQ'\over \cP'} \ = \ 
{1\over a_0} \ {{}\atop -}
{2b_1\over a_1} \ {{}\atop -} 
{b_2\over a_2} \ {{}\atop - \cdots} 
$$
which is a $ q $-analogue of Theorem B in Section 1.\qquad\endproof

\heading{5}{Ordinary cases  $ s=q^3, q^4 , \ldots $}
By substituting  $ s = q^m , m $  integer greater then 2, into
(3.5) of Theorem 4 we obtain

\thm{Corollary 9}  For  $ s = q^m $,  $ m = 3,4, \ldots $
$$
\leqalignno{
& {1\over a_0} \ {{}\atop -}
{b_1\over a_1} \ {{}\atop -} 
{b_2\over a_2} \ {{}\atop - \cdots}  \ = \ 
{ q^{(m-3)/2} (1 - q^{m-1}) ( 1 - {a\over q}) \over
( 1 - {q\over a}^{\scriptscriptstyle m-1})
( 1 - {a\over b})
( 1 - {a\over c})
( 1 - {a\over d})
( 1 - {a\over e})
( 1 - {a\over f})}& (5.1)\cr
\times \ & 
\bigg[ \phi \Big(
 {q\over a}; {q\over b}, {q\over c}, {q\over d}, {q\over e},
{q\over f}, q^{-m+2} , q \Big) -
{\tW_2\over \tW_1} \ 
 { (a)_\infty (aq)_\infty  (q)_\infty \over
({q^m\over a})_\infty
({q^{m-1}\over a})_\infty  
(q^m)_\infty} \ 
{1\over (1-q^{m-1})} \cr
& \hskip 1in \times \
{({b\over a} q^{m-1} , \ 
{c\over a} q^{m-1} , \ 
{d\over a} q^{m-1} , \ 
{e\over a} q^{m-1} , \ 
{f\over a} q^{m-1})_\infty\over
( {bq\over a},\  {cq\over a},\  {dq\over a}, \ {eq\over a}, \ 
{fq\over a})_\infty } \bigg] \cr}
$$
where  $ a_n $,  $ b_n $  are given by (2.7) and (2.8) for
$ s = q^m $,  $ m \ge 3 $.\endthm

Note that if we substitute  $ s = q $$(m=1) $  into (5.1) the right side
reduces to
$$
{q\over a^2} \ {\tW_2\over \tW_1} \ 
{ ({a\over q})_\infty (aq)_\infty \over 
({1\over a})_\infty ({1\over a})_\infty }
\leqno(5.2)
$$
and (5.1) agrees with (4.5) when the indeterminacy in
$ a_0 $  and  $ b_1 $  is taken into account, that is
$$
\eqalign{
\lim\limits_{s\to q} \ a_0 (s) \ = \ & \lim\limits_{n\to 0} \ 
a_n ( s = q) \cr
\lim\limits_{s\to q} \ b_1 (s) \ = \ & 2\lim\limits_{n\to 1} \ 
b_n ( s = q) \ .\cr}
\leqno(5.3)
$$
It is the  $ a_0 $  and  $ b_1 $  on the left of (5.3) that 
should occur in (5.1) with  $ s = q $
while it is  $ \lim\limits_{n\to 0} \ a_n (s=q) $  and
$ \lim\limits_{n\to 1} \ b_n (s=q) $  which are the  $ a_0 $
and  $ b_1 $  that occur in (4.5).

Similarly for  $ s = q^2 $,  the right side of (5.1) becomes
$$
- {a\over q^{3\over 2}} \ 
{(1-q)\over (1-{a\over b}) 
 (1-{a\over c}) 
 (1-{a\over d}) 
 (1-{a\over e}) 
 (1-{a\over f}) } \ 
 \left[ 1 -
 { (a)_\infty (aq)_\infty  \over
 ({q^2\over a})_\infty ({q\over a})_\infty } \ 
 {\tW_2\over \tW_1 }\right] \ .
 \leqno(5.4)
 $$
This agrees with (4.1) since
$$
\eqalignno{
\lim_{s\to q^2} \ a_0(s) \ = \ &
\lim_{n\to 0} \ a_n (s=q^2) +
{aq^{1\over 2} (1-{b\over a}) (1-{c\over a}) (1-{d\over a})
 (1-{e\over a}) (1-{f\over a})\over 2 (1-q) } \cr
\lim_{s\to q^2} \ b_1(s) \ = \ &
\lim_{n\to 1} \ b_n (s=q^2)\ . \cr}
$$

One can also consider the terminating case of (5.1) by taking
into account Lemma 6.
\smallskip

{\it Remarks} 
\meti{1) } If in Section 2 we make the replacements 
$ a \to \lambda a $, 
$ b \to {\lambda q\over b} $, 
$ c \to {\lambda q\over c} $, 
$ d \to {\lambda q\over d} $, 
$ e \to  ae^{i\theta} $,
$ f \to  ae^{-i\theta} $
and let  $ \lambda \to \infty $,  then we obtain solutions to
the recurrence for Askey-Wilson polynomials [1] with  
$ s = abcd = q^m $,  $ m = 1,2, \ldots $ [7].  By applying the 
above limit to Corollary 7,8, and 9 we recover equations 
(22), (23) and (24) respectively of Gupta and Masson [7].  Note that
[7, (22), (23)] give the  $ q $-analogue of Ramanujan's Entries
35 and 39 [3], [13], while Corollaries 7 and 8 are the
$ q $-analogues of Corollaries 6 and 7 of [15].

\meti{2) } The Corollary 8 case  $ s = q $  is particularly
interesting since the approximants of the continued fraction
$$
 {1\over a_0} \ {{}\atop -}
{2b_1\over a_1} \ {{}\atop -} 
{b_2\over a_2} \ {{}\atop - \cdots}  \ = \ 
{q (aq)_\infty ({a\over q})_\infty \over
a^2 ({1\over a})_\infty ({1\over a})_\infty } \ 
{\tW_2\over \tW_1}
$$
are then given explicitly in terms of  $ X_n^{(1)} $  and
$ X_n^{(2)} $.  To see this we note that the initial conditions
$$
\eqalign{
2X_1^{(1)} - a_0 X_0^{(1)} \ = \ & 0 \cr
\ X_2^{(2)} - a_1 X_1^{(2)} \ = \ & 0 \cr}
$$
which follow from (2.7), (2.8), (2.9) and (2.12) when  $ s = q $,
together with the recurrence (2.6), imply that for  $ s= q$
$$
 {1\over a_0} \ {{}\atop -}
{2b_1\over a_1} \ {{}\atop -} 
{b_2\over a_2} \ {{}\atop - \cdots - }  {b_n\over a_n } \ = \ 
{X_{n+1}^{(2)} X_0^{(1)} \over
2X_{n+1}^{(1)} X_1^{(2)} } \ , \qquad n \ge 0 \ .
$$
\meti{3) } In the limit as  $ m \to \infty $  $ (s = q^m \to 0) $
Corollary 9 yields a new continued fraction result given by
$$
 \eqalignno{
 & {1\over c_0} \ {{}\atop -}
{d_1\over c_1} \ {{}\atop -} 
{d_2\over c_2} \ {{}\atop - \cdots}  \ = \ 
{(1-{a\over q})\over
q(1-{a\over b}) (1-{a\over c}) (1-{a\over d}) ( 1-{a\over e}) } \cr
& \qquad \times
\Big[ W \big( {q\over a}; \ 
{q\over b}, \ {q\over c} , \ {q\over d} , \ {q\over e}, q \big) -
R \Big] \ , \cr}
$$
$$
\eqalignno{
R \  = \  &{\big(q, a, {q^2\over a}, \ 
{de\over a} , \  
{dc\over a}, \ 
{ec\over a}\big)_\infty \over
({dq\over a}, \ {eq\over a}, \  {cq\over a}, \ {dec\over aq}, \  
{aq\over b}, \ b\big)_\infty} 
\times  \bigg\{
{}_3\phi_2 
\Big( {{{q\over d}, {q\over e}, {q\over c}}\atop
{q {b\over a}, {q^2 a\over cde}}}\ ; \ q,b \Big)  \cr
& + {({q\over d}, {q\over e}, {q\over c} , {bcde\over a^2} , {aq\over b},
{b\over a}, {cde\over a^2} , {a^2q\over cde}, {dec\over aq} )_\infty
\over
( {ec\over a}, {qb\over a}, {qa\over cde }, {q\over a }, a, 
{bcde\over qa^2 }, {a^2 q^2\over bcde}, {de\over a }, {dc\over a} )
_\infty }\cr
& \qquad \times  {}_3\phi_2 
\Big({{ {de\over a}, {dc\over a}, {ec\over a}} \atop
{{bcde\over a^2}, {dec\over a}}}\, ; \ q, b\Big) \bigg\} \Bigg/
{}_3\phi_2 
\Big( {{{aq\over bc}, {aq\over bd}, {aq\over be}} \atop
{{aq\over b}, {a^2q^2\over bcde}}} \, ; \ q,b \Big) \cr}
$$
where
$$
\leqalignno{
c_n \ = \ & \Big[
-\big( 1 - {a\over b} q^{n+1} \big)
\big( 1 - {a\over c} q^{n+1} \big)
\big( 1 - {a\over d} q^{n+1} \big)
\big( 1 - {a\over e} q^{n+1} \big) &(5.6)\cr
& - q (1-q^n)  
 (1-aq^n) 
 (1-aq^{n+1}) 
\big (1-{a^2 q^{n+1}\over bcde}\big)  \cr
& + a^2 q^{2n+2} 
{(1-b)(1-c)(1-d)(1-e)\over bcde} \Big] \Big/ (1-aq^{n+1}) \ , \cr
d_n \ = \ &
q(1-q^n) 
\big(1 -{a\over b} q^n\big)
\big(1 -{a\over c} q^n\big)
\big(1 -{a\over d} q^n\big)
\big(1 -{a\over e} q^n\big)
\big(1-{a^2 q^{n+1}\over bcde}\big) \ . \cr}
$$
A direct proof is obtained by applying our methods to the
contiguous relation
$$
\leqalignno{
& q\big(1-{1\over f}\big)
\big(1-{a^2 q\over bcdef}\big)
\big(1-{a\over f}\big)
\big(1-{aq\over f}\big)
[ W (f +) - W] & (5.7)\cr
+ \ & 
\big(1-{aq\over fb}\big)
\big(1-{aq\over fc}\big)
\big(1-{aq\over fd}\big)
\big(1-{aq\over fe}\big)
[ W (f -) - W] \cr
+ \ &
{a^2 q^2\over bcdef^2} \ 
(1-b) (1-c) (1-d) (1-e) W \ = \ 0 \cr}
$$
where
$$
\leqalignno{
W \ = \ & W ( a; \ b,c,d,e,f) & (5.8) \cr
= \ & \ephis \Big ({{a, q\sqrt a, - q\sqrt a, b, c, d, e, f}
\atop { \sqrt a, - \sqrt a, {aq\over b} , {aq\over c}, {aq\over d},
{aq\over e}, {aq\over f}}} \ ; \ 
q, {a^2 q^2\over bcdef} \Big) \ . \cr}
$$
The contiguous relation (5.7) is obtained from (2.5) by taking the 
limit  $ g\to 0 $  with  $ fg = {a^3 q^2\over bcdeh}  $  and then
replacing  $ h = q^{-n} $  by  $ f $.  However the termination
of the above  $ \ephis $  is not necessary.

\Refs

\ref 1\\
{\smc R. Askey and J. Wilson},
{\it Some basic hypergeometric orthogonal polynomials that generalize
Jacobi polynomials},
Memoirs Amer. Math. Soc.,  319 (1985), pp.~1--55.
\endref

\ref 2\\
{\smc W.N. Bailey},
{\it Generalized Hypergeometric Series,}
Cambridge Univ. Press, London, 1935.
\endref

\ref 3\\
{\smc B.C. Berndt, R.L. Lamphere and B.M. Wilson},
{\it Chapter 12 of Ramanujan's second note-book:  
continued fractions},
Rocky Mtn. J. of Math.,  15 (1985), pp.~235--310.
\endref

\ref 4\\
{\smc G. Gasper and M. Rahman},
{\it Basic Hypergeometric Series,}
Cambridge Univ. Press, Cambridge, 1990.
\endref

\ref 5\\
{\smc D.P. Gupta, M.E.H. Ismail and D.R. Masson},
{\it Associated continuous Hahn polynomials},
Canad. J. of Math., 43 (1991), pp.~1263--1280.
\endref

\ref 6\\
{\smc D.P. Gupta, M.E.H. Ismail and D.R. Masson},
{\it Contiguous relations, Basic Hypergeometric functions and
orthogonal polynomials II, Associated Big $ q $-Jacobi polynomials},
J. of Math. Anal. and Applications, 171 (1992), pp.~477--497.
\endref

\ref 7\\
{\smc D.P. Gupta and D.R. Masson},
{\it Exceptional  $ q $-Askey-Wilson polynomials and 
continued fractions},
Proc. A.M.S., 112 (1991), pp.~717--727.
\endref

\ref 8\\
{\smc M.E.H. Ismail, J. Letessier, G. Valent and J. Wimp},
{\it Two families of associated Wilson polynomials},
Can. J. Math.,  42 (1990), pp.~659--695.
\endref

\ref 9\\
{\smc M.E.H. Ismail and M. Rahman},
{\it Associated Askey-Wilson polynomials},
Trans. Amer. Math. Soc.,  328 (1991), pp.~201--239.
\endref

\ref 10\\
{\smc L. Jacobsen}, 
{\it Domains of validity for some of Ramanujan's continued fraction 
formulas},
J. Math. Anal. and Applications, 143 (1989), pp.~412--437.
\endref

\ref 11\\
{\smc W.B. Jones and W.J. Thron}, {\it Continued Fractions: Analytic
Theory and Applications,} 
Addison-Wesley, Reading, Mass., 1980.
\endref

\ref 12\\
{\smc D.R. Masson}, 
{\it Some continued fractions of Ramanujan and Meixner-Pollaczck
polynomials},
Canad. Math. Bull.,  32 (1989), pp.~177--181.
\endref

\ref 13\\
{\smc D.R. Masson}, 
{\it Wilson polynomials and some continued fractions
of Ramanujan},
Rocky Mountain J. of Math.,  21 (1991), pp.~489--499.
\endref

\ref 14\\
{\smc D.R. Masson}, 
{\it Associated Wilson polynomials}, 
Constructive Approximation,  7 (1991), pp.~521--534.
\endref

\ref 15\\
{\smc D.R. Masson}, 
{\it A generalization of Ramanujan's best theorem on continued
fractions},
Can. Math. Reports of the Acad. of Sci., 13 (1991), pp.~167--172.
\endref

\ref 16\\
{\smc S. Ramanujan},
{\it Notebooks}, Tata Institute, Bombay, 1957.
\endref

\ref 17\\
{\smc G.N. Watson},
{\it Ramanujan's continued fraction},
Proc. Cam. Philos. Soc., 31 (1935), pp.~7--17.
\endref

\ref 18\\
{\smc J.A. Wilson},
{\it Hypergeometric series, recurrence relations and some new 
orthogonal polynomials},
Ph.D. thesis, University of Wisconsin-Madison, 1978.
\endref

\ref 19\\
{\smc L.C. Zhang},
{\it Ramanujan's continued fractions for products of gamma functions}, 
to appear.
\endref
\bye